\documentclass[12pt]{article}
\usepackage{amssymb}
\begin{document}

\renewcommand{\thefootnote}{\arabic{footnote}}
\newcommand{\blambda}{\mbox{\boldmath {$\lambda$}}}
\newcommand{\bgamma}{\mbox{\boldmath {$\gamma$}}}
\newtheorem{thm}{Theorem}[section]
\newtheorem{corollary}{Corollary}[section]
\newtheorem{conjecture}{Conjecture}[section]
\newtheorem{prop}{Proposition}[section]
\newtheorem{lemma}{Lemma}[section]
\newtheorem{definition}{Definition}[section]
\newtheorem{example}{Example}[section]

\def\blackslug{\hbox{\hskip 1pt \vrule width 4pt height 8pt depth 1.5pt
\hskip 1pt}}
\def\proof{\par\penalty-1000\vskip .5 pt\noindent{\bf Proof\/: }}
\def\Reals{{\bf R}}
\def\vect{\bf}
\def\per{\bf Per}

\begin{center}{\Large \bf Multivariate normal approximations by Stein's method and size bias couplings}\end{center}

\vspace{3mm}

\begin{center}{Larry Goldstein$^*$  and Yosef Rinott$^{\dag}$}\end{center}

\vspace{3mm}

\vspace{2.5in}

\abstract  Stein's method is used to obtain two theorems on multivariate normal approximation. Our main theorem, Theorem \ref{main}, provides a bound on the distance to normality for any nonnegative random vector.  Theorem \ref{main} requires multivariate size bias coupling, which we discuss in studying the approximation
of distributions of sums of dependent random vectors. In the univariate case,
we briefly illustrate this approach for certain sums of nonlinear functions of multivariate normal variables. As a second illustration, we show
that the multivariate distribution counting the number of vertices with given degrees in certain random graphs is asymptotically multivariate normal and obtain a bound on the rate of convergence. Both examples demonstrate that
this approach may be suitable for situations involving non-local dependence.
We also present Theorem \ref{cor2} for sums of vectors having a local type of dependence. We apply this theorem to obtain a multivariate normal
approximation for the distribution of the random $p$-vector
which counts the number of edges in a fixed graph both of whose vertices
have the same given color when each vertex is colored by one
of $p$ colors independently. All normal approximation results
presented here do not require an ordering of
the summands related to the dependence structure. This is in contrast to
hypotheses of classical central limit theorems and examples, which involve e.g., martingale, Markov chain, or various mixing assumptions.

\endabstract

\vfill

\noindent \underline{\hspace{2.5in}}



\noindent $^*${\footnotesize Department of Mathematics DRB-155, USC, Los Angeles, CA 90089-1113. This work was supported in part by NSF grant DMS 90-05833}

\noindent $^{\dag}${\footnotesize Department of Mathematics, UCSD, La Jolla, CA  92093. This work was supported in part by NSF grant DMS 92-05759.}

\newpage

\section{Introduction}
Stein's method has been successful in assessing the quality of normal and Poisson approximations under various dependence structures.
See Stein (1972), Stein (1986), Barbour, Holst and Janson (1992), and
references therein. Significant multivariate (or functional)
versions of Stein's
method appear for example in Barbour (1990), and
G\"otze (1991). An important part
of Stein's method is the construction of auxiliary random variables (coupling), which
are used in the computation of the bounds
on the distance between a given random variable and its normal or Poisson
approximant.

The coupling variables constructed in the application of Stein's method may not appear explicitly in the final bounds.
Applications of results where the auxiliary random variables do appear in the bounds require  their explicit construction. Although this feature may make
such results more difficult to use, couplings yielding useful
bounds may often be found where other methods seem to fail. For the Poisson,
theorems of this nature may be found in Barbour, Holst and Janson (1992) and references within.
In the main result of this paper, Theorem \ref{main},  we obtain bounds for normal approximations in terms of such couplings, and provide general guidelines and methods for their construction so that these methods may be applied. The
couplings studied here, an instance of a construction of a joint distribution with given
marginals, is of independent interest. It is known that constructions of
multivariate couplings may be problematic, see, e.g., Dall'Aglio, Kotz and
Salinetti (1991). Nevertheless we are able to provide methods for the required multidimensional coupling constructions, which we illustrate in two applications involving nonlocal dependence, Theorems \ref{gauss} and \ref{ex1}.

By the same techniques used to prove Theorem \ref{main}, Stein's method and the analysis of the properties of the solution to the partial differential equation (\ref{DE}), we obtain Theorem \ref{cor2}, a result complementary to Theorem \ref{main}. Theorem \ref{cor2} provides a multivariate normal approximation under conditions of local dependence. Unlike Theorem \ref{main}, coupling variables do not appear explicitly in Theorem \ref{cor2}.

In order to introduce the couplings needed for the proof and applications of Theorem \ref{main}, we require the following definition.
\begin{definition}
\label{defodsb}
Given a nonnegative random variable $W$ with distribution $dF(w)$ and mean $\lambda$, $W^*$ is said to have the $W$-size biased distribution if it has distribution $wdF(w)/\lambda$.
\end{definition}
Note that the distribution of $W^*$ may be characterized by the relation
\begin{equation}
\label{onedsb}
EWG(W)=\lambda EG(W^*)
\end{equation}
for all functions $G$ for which the expectations exist.
Size biased distributions are well known in sampling theory and renewal theory, for example.
The following one dimensional version of our main result illustrates the relevance of size biased
coupling to normal approximations.
\begin{thm}
\label{unimain}
Let $W$ be a nonnegative random variable with mean $EW=\lambda$, variance $\sigma^2=\mbox{Var}(W)$, and let $W^*$ be jointly
distributed with $W$, having the $W$-size biased distribution.
Then for any piecewise continuously differentiable $h$,
\begin{eqnarray}
\lefteqn{|Eh(\frac{W-\lambda}{\sigma})-\Phi h| } \nonumber \\
&& \leq 2||h||\frac{\lambda}{\sigma^2}
\sqrt{VarE(W^*-W|W)}+||h'||\frac{\lambda}{\sigma^3}
E(W^*-W)^2,
\end{eqnarray}
where $||\cdot||$ denotes the supremum norm, and $\Phi h=Eh(Z)$ with $Z$ a standard normal variate.
\end{thm}

Theorem \ref{unimain} is an extension of a
result of Baldi, Rinott and Stein (1989).
The theorem requires the construction of $W^*$ on a joint space with $W$; hence, obtaining good bounds in any particular application depends on the construction of a $W^*$ which will be
close to $W$ in an appropriate sense. However, since the resulting bound is valid for any construction for which the marginal distribution of $W^*$
coincides with the $W$-size biased distribution, one has the flexibility to choose constructions which result in good computable bounds.

Here is a brief description of a method which leads
to a construction of a $W$-size biased variate $W^*$ when $W=X_1+\cdots+X_n$ is a sum of random variables.
To begin,  if $X_1,\ldots,X_n$ are iid nonnegative random variables with finite mean, then $W^*$ can be constructed by replacing
any single summand, say, $X_1$ by an independent variable $X_1^*$ with the $X_1$-size biased
distribution, i.e. $W^*=X_1^*+X_2+\cdots+X_n$. More generally, if $W$ is a sum of non-iid variates, then a like construction of $W^*$ may be given by replacing $X_I$ by $X_I^*$, where the random index
is chosen independently with $P(I=i)=EX_i/\sum EX_j$, and adjusting the remaining variables
to their conditional distribution given the new value of $X_I$.

A special case of this idea is Midzuno's procedure (e.g. Cochran (1977)), where a size biased variable is used to obtain unbiased ratio estimators in finite population sampling.
To describe Midzuno's procedure, let nonnegative ``sizes'' $X_1,\ldots,X_n$ be obtained by sampling from a finite population without replacement, and $W$ be their sum. Then, $W^*$ is realized by sampling the first variate in proportion to its size, removing it from the population, and sampling the other variables
without replacement from the population that remains, that is, sampling from the resulting conditional distribution.

Further flexibility is obtained by realizing that for any representation of $W$ in the form
\begin{equation}
\label{fun}
W=\psi_1(U_1)+\cdots+\psi_n(U_n),
\end{equation}
the above construction of $W^*$ may be accomplished by choosing a random index $I$
such that $P(I=i)=E\psi_i(U_i)/\sum_{j=1}^nE\psi_j(U_j)$, and if $I=i$, replacing $U_i$ by an independent variable with distribution $\psi_i(u)P(U_i \in du)/E\psi_i(U_i)$, and adjusting the remaining $U$ variables. Therefore, the theorem may
be applied whenever one can find a transformation such that the variables $U_1,\ldots,U_n$ have a dependence
structure that allows the computation of the conditional distribution required in the bounds.  Further details and examples of these size bias coupling constructions and their applications will be provided
in Sections \ref{secsize} and \ref{ex}.

For the multivariate case, we need a more general notion of size
biasing, and we replace the
$^*$ notation by a superscript ${\beta}$ in order to identify
in which ``coordinate" or variable the variates are size biased.

\begin{definition}
\label{size}
Let ${\cal I}$ an  arbitrary index set and let ${\bf X}=\{ X_{\alpha}: \alpha \in {\cal I}\}$ be a collection of nonnegative random variables with joint
distribution $dF({\bf x})$ and means $EX_{\alpha}=\lambda_{\alpha}$. For $\beta \in {\cal I}$,
we say that ${\bf X}^{\beta}= \{ X_{\alpha}^{\beta}: \alpha \in {\cal I}\}$ has
the ${\bf X}$-size biased distribution in the $\beta^{\rm th}$ coordinate if ${\bf X}^{\beta}$ has the joint distribution $x_{\beta}dF({\bf x})/\lambda_{\beta}.$
\end{definition}

\noindent The distribution of ${\bf X}^{\beta}$ is characterized by the relations
\begin{equation}
\label{genrest}
E X_{\beta} G({\bf X}) = \lambda_{\beta} EG({\bf X}^{\beta})
\end{equation}
for all functions $G$ for which the above expectations exist.
When the function $G$ depends on ${\bf X}$ only through $X_{\beta}$,
equation (\ref{genrest}) yields $EX_{\beta}G(X_{\beta})=\lambda_{\beta}EG(X_{\beta}^{\beta})$,
hence, comparing to (\ref{onedsb}), we see that the $\beta^{\rm th}$ coordinate of ${\bf X}^{\beta}$, that is, the variate $X_{\beta}^{\beta}$,
has the $X_{\beta}$-size biased distribution in the sense of definition \ref{defodsb}.

By considering the case where the collection ${\bf X}$ consists of only
a single random variable, we see that equation (\ref{genrest}) reduces to
(\ref{onedsb}); hence definition \ref{defodsb} is a special case of defintion \ref{size}.

We will apply Definition \ref{size} to a vector
${\bf W}=(W_1,\ldots,W_p)\in R^p $ by identifying it with the collection $\{ W_j: j \in {\cal I}\}$ with ${\cal I}=\{1,\ldots,p\}$. Letting $EW_i=\lambda_i$, we see that the vector ${\bf W}^i=(W_1^i,\ldots,W_p^i)$
is characterized by
\begin{equation}
\label{rest}
EW_iG({\bf W})=\lambda_i EG({\bf W}^i).
\end{equation}

Relation (\ref{rest}) leads to a multivariate normal approximation
theorem, for which we introduce the following notation (see e.g.,
Horn and Johnson  (1985)). Given a vector ${\bf a}$ in $R^p$, let
$||{\bf a}||=\max_{1 \leq i \leq p}|a_i|$. Given a $p \times p$
matrix $A=(a_{ij})$ we set $||A||=\max_{1\leq i,j \leq p} \,
|a_{ij}|$\,, and more generally for any array, $||\cdot||$ will
denote its maximal absolute value. For an array of functions, say
$A({\bf w})=\{a_{i}({\bf w})\}$, where $i$ could stand for a
multiple index, $||A||=\sup_{\bf w}\max_i|a_{i}({\bf w})|$. For a
smooth function $h:R^p \rightarrow R$ we let $\nabla h$ or $Dh$
denote the vector of first partial derivatives of $h$, $D^2h$ the
usual Hessian matrix of second order partial derivatives and
$D^kh$ the $k^{\rm th}$ derivative of $h$ in general.

Our main multivariate result is the following theorem:
\begin{thm}
\label{main}
Let ${\bf W}$ be a random vector in $R^p$ with nonnegative components. Set
$\blambda = (\lambda_1, \ldots, \lambda_p) = E{\bf W}$, and assume $\mbox{Var}{\bf W} = \Sigma = (\sigma_{ij})\,$ is invertible.
For each $i=1,\ldots,p$ let $({\bf W}, {\bf W}^i)$ be
random vectors defined on a joint probability space with ${\bf W}^i$ having the ${\bf W}$-
size biased distribution in the $i^{\rm th}$ coordinate as in (\ref{rest}).
Let $h : R^p \rightarrow R$ be a function having bounded mixed partial
derivatives up to order 3.
Let $\Phi h = Eh({\bf Z})$, where ${\bf Z}$ denotes a standard (mean zero, covariance $I$)
normal vector in $R^p$.
Then
\begin{eqnarray}
\label{bounda}
\lefteqn{|Eh(\Sigma^{-1/2}({\bf W} - {\blambda})) - \Phi h|  \leq}
 \nonumber \\
&&\frac{p^2}{2}||\Sigma^{-1/2}||^2||D^2h||
\sum_{i=1}^p \sum_{j=1}^p\lambda_i  \sqrt{\mbox{Var}E[\,W_j^i-W_j\, |\, {\bf W}]}\nonumber\\
&& + \frac{1}{2}\frac{p^3}{3}||\Sigma^{-1/2}||^3||D^3h||
\sum_{i=1}^p \sum_{j=1}^p\sum_{k=1}^p
\lambda_iE|(W_j^i-W_j)(W_k^i-W_k)|.
\end{eqnarray}
\end{thm}

Note that the theorem does not require the joint construction of $({\bf W}^1,\ldots,{\bf W}^p)$.
Although Theorems \ref{unimain} and \ref{main} are stated for nonnegative variates, they may be applied to general variates by translation and truncation.

In Section \ref{secsize}, we discuss the construction
of the vectors ${\bf W}^i$
required for
Theorem \ref{main} when the components of ${\bf W}$
are sums of dependent random variables.  Specifically, when ${\bf X}=\{ X_{\alpha} , \alpha \in {\cal I}\}$ is any collection of nonnegative random
variables, and $A_1,\ldots,A_p$ are any subsets of ${\cal I}$, we may apply
Theorem \ref{main} to the vector ${\bf W}=(W_1,\ldots,W_p)$ where $W_j= \sum_{\alpha \in A_j}X_{\alpha}$. In particular, we obtain a result for a sum
${\bf W}=(W_1,\ldots,W_p)=\sum_{u=1}^n{\bf X}_u$ of nonnegative dependent random vectors, ${\bf X}_u = (X_{u1}, \ldots, X_{up}), \, u=1, \ldots, n,$ by letting ${\cal I}$ be a set of double indices and $A_j=\{1,\ldots,n \} \times j$.

We briefly indicate how size biased variables arise in one dimensional normal approximations.  Given a random variable $W$ and a test function $h$, one can compute $E[h(\frac{W-\lambda}{\sigma})-\Phi h]$
by computing $E[f'(W)-\frac{(W-\lambda)}{\sigma^2}f(W)]$ where $f$ is the
bounded solution of
the Stein equation
\begin{equation}
\label{charles}
f'(w)-\frac{(w-\lambda)}{\sigma^2}f(w)=h(\frac{w-\lambda}{\sigma})-\Phi h.
\end{equation}
If $W^*$ has the $W$-size biased distribution, and therefore satisfies $EWf(W) = EWEf(W^*)$,
we obtain
\[ Eh(\frac{W-\lambda}{\sigma})-\Phi h
=E[f'(W)-\frac{(W-\lambda)}{\sigma^2}f(W)]=E[f'(W)-\frac{\lambda}{\sigma^2}(f(W^*)-f(W))].\]
Taylor expansion of $E[f(W^*)-f(W)]$ is then the first step in obtaining the bound in Theorem \ref{unimain}. Note that the one dimensional versions of the theorems are not exactly special cases of their multivariate counterparts.  In the multivariate case, equation (\ref{charles}) will be replaced by the partial differential equation (\ref{DE}), resulting in different
orders of the derivatives of $h$ appearing in the one and multidimensional theorems.

The size biased coupling approach
handles cases where there is global dependence among the
summand variables. In contrast the following univariate and
multivariate results not based on size biased couplings are
very useful in cases
of local dependence.
The following theorem is due to Stein (1986); our Theorem \ref{cor2} is a multivariate version.
\begin{thm}
\label{unicor2}
Let $ X_v, \, v=1, \ldots, n, $ be random
variables with $EX_v = 0.$
Let $S_v, \, v=1,\ldots,n$ be subsets  of $\{1,\ldots, n\}$ and set $W = \sum_{v=1}^nX_v$ and denote
\[ \sum_{v=1}^n\sum_{u \in S_v}
E X_{v}X_{u}= \sigma^2,\]
assuming $\sigma^2>0$. Then for any $h : R \rightarrow R$ which is continuous and piecewise continuously differentiable,

\begin{eqnarray}
\label{unisecbounda}
\lefteqn{|Eh(W/\sigma) - \Phi h|  \leq
\frac{2||h||}{\sigma^2} \sqrt{E\left\{\sum_{v=1}^n\sum_{u \in S_v}
(X_{v}X_{u}-EX_{v}X_{u})\right\}^2}} \nonumber \\
&& +  \,\,\sqrt{\frac{\pi}{2}}\frac{ ||h||}{\sigma} \sum_{v=1}^n\,
E| E[X_{v}\,|\,X_u\,:\, u \not\in S_v]\,| \nonumber \\
&&+\frac{1}{\sigma^3} ||h'||
\sum_{v=1}^nE|X_{v}
\sum_{u \in S_v}X_{u}\sum_{t \in S_v}X_{t}|.
\end{eqnarray}
\end{thm}
In typical applications of Theorem \ref{unicor2}, we have $X_v$ independent of $\{X_u:u \not \in S_v\}$ and the second term of the bound in (\ref{unisecbounda}) vanishes. In this case we may view $S_v$ as
a {\em dependency neighborhood} of $X_v$.  Generally, the bound in (\ref{unisecbounda}) is small
if these neighborhoods are small, so this theorem is useful when the dependence is local.

The following result, which is particularly useful for normal approximations of sums of locally dependent random vectors, extends Theorem \ref{unicor2} to the multivariate case.
\begin{thm}
\label{cor2}
Let $\{ X_{\alpha} , \alpha \in {\cal I}\}$ be random
variables with $EX_{\alpha} = 0$. Let $A_1,\ldots,A_p$ be subsets of ${\cal I}$, and set ${\bf W}=(W_1,\ldots,W_p)$ where $W_j= \sum_{\alpha \in A_j}X_{\alpha}$.
For each $\alpha \in {\cal I}$
let $S_{\alpha} \subseteq {\cal I}$, and assume that $\Sigma = (\sigma_{ij})$ is symmetric positive definite, where
\[ \sigma_{ij}=\sum_{\alpha \in A_i}\sum_{\beta \in A_j \cap S_{\alpha}}
E X_{\alpha}X_{\beta}, \quad i,j = 1, \ldots, p.\]
Let $h : R^p \rightarrow R$ be a function having bounded mixed partial
derivatives up to order 3,
and $\Phi h = Eh({\bf Z})$ where ${\bf Z}$ denotes a standard normal vector in $R^p$.
Then
\begin{eqnarray}
\label{secbounda}
\lefteqn{|Eh(\Sigma^{-1/2}{\bf W}) - \Phi h|  \leq
\frac{p^2}{2}||\Sigma^{-1/2}||^2||D^2h||
\,\sum_{i=1}^p\sum_{j=1}^p\, \sqrt{E\left\{\sum_{\alpha \in A_i}\sum_{\beta \in A_j \cap S_{\alpha}}
(X_{\alpha}X_{\beta}-EX_{\alpha}X_{\beta})\right\}^2}} \nonumber \\
&& +  p||\Sigma^{-1/2}|| \,\, ||Dh||\,\sum_{i=1}^p\sum_{\alpha \in A_i}\,
E| E[X_{\alpha}\,|\,X_{\beta}\,:\, \beta \not\in S_{\alpha}]\,| \nonumber \\
&&+
\frac{1}{2}\frac{p^3}{3}||\Sigma^{-1/2}||^3||D^3h||
\,\sum_{i=1}^p\sum_{j=1}^p\sum_{k=1}^p
\sum_{\alpha \in A_i}E|X_{\alpha}
\sum_{\beta \in A_j \cap S_{\alpha}}X_{\beta}\sum_{\gamma \in A_k \cap S_{\alpha}}X_{\gamma}|.
\end{eqnarray}
\end{thm}

Note that $\sigma^2$ of Theorem \ref{unicor2} and $\Sigma$ of Theorem \ref{cor2} are not necessarily equal to the covariance of $W$ and the
covariance matrix of ${\bf W}$, respectively.
In Theorem \ref{cor2}, the symmetry of ${\bf \Sigma}$ is guaranteed if the sets
$S_{\alpha} \subseteq {\cal I}$ are symmetric in the sense that $\beta \in S_{\alpha}$ if and only if $\alpha \in S_{\beta}$. In particular, in applying Theorem \ref{cor2} to a sum of mean zero random vectors, ${\bf W}=\sum_{u=1} ^n {\bf X}_u$, where ${\bf X}_u = (X_{u1}, \ldots, X_{up})$ for $u=1, \ldots, n,$ it is natural to take neighborhoods of the form $S_{\alpha}=S_{(u,i)}=T_u \times \{1,\ldots,p\}$,
where $T_u$ are symmetric subsets of $\{1,\ldots, n\}$.

In our applications, the sets $S_{\alpha}$ contain $\{ \beta \in {\cal I}: \mbox{Cov}(X_{\beta},X_{\alpha})\not = 0\};$
in any such case, $\Sigma =\mbox{Cov}({\bf W})$. In particular, $\Sigma =\mbox{Cov}({\bf W})$ if $X_{\beta}$ is independent of $X_{\alpha}$ for every $\beta \not \in S_{\alpha}$. In the general case, the above
(somewhat unusual) choice of $\Sigma$ simplifies the form
of the bound. With the more natural choice
$\Sigma=\mbox{Cov}({\bf W})$,
the present technique applies to yield
a version of the above theorem,
but an additional term in the bound may result.

Applications of the above theorems are given in Section \ref{ex}.
As an illustration involving non-local dependence, we apply Theorem \ref{unimain} and its multivariate extension, Theorem \ref{main}, to show that the multivariate distribution counting the number of vertices with given degrees in certain random graphs is asymptotically multivariate normal and obtain a bound on the rate of convergence. To illustrate a case of local dependence, we apply Theorem \ref{unicor2}, and its multivariate extension,
Theorem \ref{cor2}, to obtain a multivariate normal
approximation for the distribution of the random $p$-vector
which counts the number of edges in a fixed graph both of whose vertices
have the same given color when each vertex is colored by one
of $p$ colors independently. Applications related to representations of $W$ as in (\ref{fun}) will be given where the $U$ variables are normal and
multinomial. The ideas and results presented here have been applied in work of
Luk (1994) in finite population sampling, and Reinert (1994) in the study of empirical measures.

The proofs of Theorems \ref{main} and \ref{cor2} are
given in Section \ref{pf}, Theorems \ref{unimain} and \ref{unicor2} may be
proved similarly.

The theorems presented here supply approximations in terms of expectations of smooth test
functions $h$, allowing our main theorems to be presented in a form where they can be readily applied under unrestrictive, simple conditions.
In the context of Stein's method, Stein (1986),  Baldi, Rinott and Stein (1989), G\"{o}tze (1991), and Rinott (1994) among others, consider also non smooth
functions, usually at the expense of added technical detail or some loss of information in the bounds. It is possible to obtain certain multivariate
version of our results for nonsmooth functions $h$ using the methodology developed in G\"{o}tze (1991), see Rinott and Rotar (1994). In the present paper, our main focus is in
the coupling structure. The issue of smooth versus non-smooth function approximation is discussed in Barbour, Karo\'{n}ski and Ruci\'{n}ski (1989).

\section{Construction of size biased couplings}
\label{secsize}
The construction of size biased variables required for the
application of Theorems \ref{unimain} and \ref{main} is the focus of this section.
While the details depend on the case at hand, this section will
provide general guidelines that extend and unify
ideas which appeared in Baldi, Rinott, Stein (1989), and Stein (1992),
where only univariate sums of zero-one variables were studied.

The following lemma is the key in the construction of
coupled variables satisfying equations (\ref{onedsb}) and (\ref{rest}) required in Theorems \ref{unimain} and \ref{main} respectively.
Readers interested only in the
univariate case may read the lemma below with ${\cal I}=A=B$.

\begin{lemma}
\label{key}
Let ${\cal I}$ be an arbitrary index set, and let ${\bf X}=\{X_{\alpha}$ :\, $\alpha \in {\cal I}\}$ be a collection of nonnegative
random variables. For any subset $B \subseteq {\cal I}$, set $X_B =
\sum_{\beta \in  B}X_{\beta}$, and
$\lambda_B=EX_B$.
Suppose $B \subseteq {\cal I}$ with $\lambda_B< \infty$,
and for $\beta \in B$  let ${\bf X}^{\beta}$ have the ${\bf X}$-size biased distribution in coordinate $\beta$ as in Definition \ref{size}.
Let ${\bf X}^B$ be a random variable distributed as the mixture of the distributions ${\bf X}^{\beta},\,\beta \in B$ with weights $\lambda_{\beta}/\lambda_B$.
Then
\begin{equation}
\label{nu}
EX_B G({\bf X})=\lambda_B EG({\bf X}^B).
\end{equation}
Hence, for any $A \subseteq {\cal I}$, if $G$ is
a function of $X_A$ only, then
\begin{equation}
\label{notasn}
EX_B G(X_A)=\lambda_B EG(X_A^B),
\end{equation}
where
\[ X_A^B= \sum_{\alpha \in A} X_{\alpha} ^B.\]
In particular, by taking $A=B$ in (\ref{notasn}) we have $EX_AG(X_A)=\lambda_A G(X_A^A)$,
and hence $X_A^A$ has the $X_A$-size biased distribution in the sense of defintion \ref{defodsb}.
\end{lemma}

{\bf Proof:} For a function $G$ on ${\bf X}$, we have $EG({\bf
X}^{\beta})=EX_{\beta}G({\bf X})/\lambda_{\beta}$ by equation
(\ref{genrest}). Multiplying by $\lambda_{\beta}/\lambda_B$ and
summing over $\beta \in B$ yields (\ref{nu}). The remainder of the
lemma now follows. $\square$

{\bf Construction of} ${\bf X}^B$:  Since ${\bf X}^B$ is a mixture of the
distributions ${\bf X}^{\beta}$ for $\beta \in B$ with weights $\lambda_{\beta}/\lambda_B$, given the collection
${\bf X}=\{X_{\alpha}$ :\, $\alpha \in {\cal I}\}$ first choose an
independent index $I \in B$
according to the distribution $P(I=\beta)=\lambda_{\beta}/\lambda_B$.
If $I=\beta$, construct $X^{\beta}_{\beta}$ to have the $X_{\beta}$-size biased distribution  $x_{\beta} P(X_{\beta} \in dx_{\beta})/\lambda_{\beta}$. If $X^{\beta}_{\beta}=x$ then
the remaining variates $X_{\alpha}^{\beta},\,\alpha \not = \beta$ are constructed so that $P({\bf X}^{\beta} \in d{\bf x})=P({\bf X} \in d{\bf x}\,| \,X_{\beta} = x)$. This construction yields
\[ P({\bf X}^\beta \in d{\bf x}) =
P({\bf X} \in d{\bf x}\,| \,X_\beta = x_\beta)\,x_\beta
P(X_\beta \in dx_\beta)/\lambda_\beta, \]
that is, ${\bf X}^\beta \sim x_{\beta}dF({\bf x})/\lambda_{\beta},$
and indeed ${\bf X}^\beta$ has the ${\bf X}$-size biased
distribution in the ${\beta}^{\rm th}$ coordinate as given in  Definition \ref{size}.

In the univariate case, with $W=\sum_{\alpha \in {\cal I}}X_\alpha=X_{\cal I}$,  and $A=B={\cal I}$, equation (\ref{notasn}) in Lemma \ref{key} shows that
a construction of $W^*$ satisfying (\ref{onedsb}) may be obtained by setting  $W^*=X_{\cal I}^{\cal I}.$ Hence $W^*$ may be constructed as follows:
a summand $X_{\beta}$ of $W$, chosen with probability
$EX_\beta/EW$, is replaced by a new value
from its size biased distribution, and the remaining summands
are adjusted to have the conditional distribution of ${\bf X}$
conditioned on the event that for the chosen $\beta$, $X_{\beta}$ takes the new value.

If the variates $\{X_{\alpha}: \alpha \in {\cal I}\}$ are
independent the last step
is not needed since by independence the conditioning is irrelevant.
In this case, the construction of a $W$-size biased variable $W^*$ reduces
to size biasing a single randomly chosen summand $X_I$.
In the case that $\{X_{\alpha} :\, \alpha \in {\cal I}\}$ are all
zero-one variates we simply have $X_\beta^\beta \equiv 1$, so
in the case that $W$ is a sum of independent zero-one random variables,
the coupling is accomplished by choosing
an index $I=\beta$ with probabilities proportional to $P(X_{\beta}=1)$, setting $X_{\beta}=1$, and leaving the remaining variates unchanged.

In the multivariate case, the connection between Lemma \ref{key} and Theorem \ref{main}
for approximating sums of random variables is as follows. Given $\{X_{\alpha} :\, \alpha \in {\cal I}\}$, let $A_1,\ldots,A_p$ be
subsets of ${\cal I}$, and set

\[ {\bf W}=(W_1,\ldots,W_p), \quad \mbox{where} \quad W_j=X_{A_j} \quad \mbox{and} \quad   {\bf W}^i=(W_1^i,\ldots,W_p^i)\quad \mbox{where} \quad W_j^i=X_{A_j}^{A_i}.\]

When $B=A_i$ and $G({\bf X})$ is a function depending on ${\bf X}$ only through ${\bf W}$, equation (\ref{nu}) yields equation (\ref{rest}). Therefore, one may
obtain the vector ${\bf W}^i$ satisfying (\ref{rest}) by constructing ${\bf X}^{A_i}$ using a random index $I \in A_i$ with $P(I=\beta)=\lambda_{\beta}/\lambda_{A_i}$ as described above.

In particular, the sum of random vectors ${\bf W}=\sum_{u=1}^n {\bf X}_u$ where ${\bf X}_u=(X_{u1},\ldots,X_{up})$, corresponds to the choice ${\cal I}=\{1,\ldots,n\} \times  \{1,\ldots,p\}$ and $A_j=\{(u,j): u=1,\ldots,n\}$. In Section \ref{ex}, we apply this multivariate construction in the setup where $X_{uj}$ is the indicator of
the event that the degree of the vertex $u$ in the random graph $K=K_{n,\pi}$ equals
a prescribed number $d_j$. Hence, $W_j$ is the number of vertices of $K$ of degree $d_j$.
The coupling in this case is accomplished as follows. Since $EX_{ui}$ are equal for indices in $A_i$, to construct ${\bf X}^{A_i}$ it is required to choose an index, say $v$, uniformly over $\{1,\ldots,n\}$ and size
bias $X_{vi}$ for this $v$. As $X_{vi}$ is an indicator, size biasing is accomplished by
replacing $X_{vi}$ by the constant 1. The above construction now requires that the remaining
variables have their original distribution conditioned on $X_{vi}=1$. If $X_{vi}$ was
initially 1, that is, if the degree of vertex $v$ was $d_i$, no change is required. Otherwise, by adding or removing randomly chosen edges as appropriate, the degree of $v$
is made to be $d_i$, thereby size biasing the indicator $X_{vi}$. This procedure results
in a new graph $K^i$, in which the other variables now have the proper conditional distribution.

The following comments pertain to the random choice of index that
appears in the above constructions. In certain cases the size
biased distribution can be constructed with a deterministic index,
however, such constructions may lead to larger bounds than those
obtained using randomization.

We specialize to $p=1$ and the expression
$\mbox{Var}E[W^*-W\,|\,W\,]$ which appears in the first term
in the bound in (\ref{bounda}). If $X_1, \ldots, X_n$ are exchangeable then it is easy to see
that $W$ may be size biased by size biasing $X_1$; that is, in the above
description one can set $I \equiv 1$ deterministically and the rest is done as above.
This is true for the example
where $W$ counts the number of vertices of degree $d$ in the random graph $K$.
Often in order to make calculations tractable it
is necessary to condition on a larger $\sigma$-field ${\cal F} \supset \sigma\{W\}$ and replace
$\mbox{Var}E[W^*-W\,|\,W\,]$ by the larger quantity
$\mbox{Var}E[W^*-W\,|\,{\cal F}\,]$. However, $\mbox{Var}E[W^*-W\,|\,{\cal F}\,]$ may
not give rise to a useful bound unless $I$ is randomized. This is the case in the random graph example when $\mbox{Var}E[W^*-W\,|\,W\,]$ is replaced by  $\mbox{Var}E[W^*-W\,|\,K\,]$.

This difficulty can be seen even in the
case of independent, identically distributed zero-one random variables, where conditioning
on ${\bf X}$ is the analog to conditioning on $K$. There, size biasing
with a random $I$ leads to $\mbox{Var}E[W^*-W\,|\,{\bf X}\,]=\mbox{Var}E[1-X_I\,|\,{\bf X}\,]=\mbox{Var}E[X_I\,|\,{\bf X}]=\mbox{Var}(W/n)$ of order
$1/n$, but setting $I \equiv 1$, the quantity
$\mbox{Var}E[W^*-W\,|\,{\bf X}\,]=\mbox{Var}(X_1),$ a constant.

The following lemma of Dembo and Rinott (1994) shows how to size bias a sum
of the form $W=\sum_{j=1} ^n \psi_j(U_j)$ by working with the argument
distribution $U_1,\ldots,U_n.$

\begin{lemma}
\label{biasfun}
Let ${\bf U} = (U_1, \ldots U_n)$ be a random $n$ vector,
and let
$\psi_i$ be nonnegative functions such that $E\psi_i(U_i) < \infty, \, i= 1, \ldots, n$.
Let ${\bf Y}^{(i)}=(Y^{(i)}_1, \ldots,Y^{(i)}_n)$ satisfy
$P({\bf Y}^{(i)} \in d{\bf y} )=
P({\bf U} \in d{\bf y})\psi_i(y_i)/E\psi_i(U_i)$.
Let $I$ be a random variable taking values in \{1, \ldots, n\},
distributed independently of all the above variables, with $P(I=i)=E\psi_i(U_i)/\sum_{j=1}^nE\psi_j(U_j)$.

Let $W=\sum_{j=1}^n \psi_j(U_j)$ have the distribution $F$. Then
$W^* = \sum_{j=1}^n \psi_j(Y^{(I)}_j)$
has the distribution $wdF(w)/\lambda$, where $\lambda=EW$.
\end{lemma}

Note that
with $F_j$ denoting the marginal distribution function of $U_j$,
the distribution of ${\bf Y}^{(i)}$ is obtained by letting
$Y^{(i)}_i$ have the marginal distribution
$\psi_i(\cdot)dF_i(\cdot)/E\psi_i(U_i)$, and if $Y^{(i)}_i=u$, letting
$(Y_1^{(i)}, \ldots, Y_{i-1}^{(i)},Y_{i+1}^{(i)}, \ldots,Y_n^{(i)})$
have the distribution of
$(U_1, \ldots, U_{i-1},U_{i+1}, \ldots,U_n)$
conditioned on $U_i=u$.

To summarize, given $W=\sum_{j=1}^n \psi_j(U_j)$,
this suggests the following:\\
{\bf Construction of} $W^*$.
Choose a random index $I$ as in the lemma. If $I=i$, let $Y_i \sim \psi_i(u)dF_i(u)/E\psi_i(U_i)$.
If $Y_i$ is assigned the value $u$, let
$(Y_1, \ldots, Y_{i-1},Y_{i+1}, \ldots,Y_n)$
have the conditional distribution of
$(U_1, \ldots, U_{i-1},U_{i+1}, \ldots,U_n)$ given $U_i=u$.
Now set $W^*  = \sum_{j=1}^n\psi_j(Y_j)$.

If $(U_1, \ldots U_n)$ are Gaussian or multinomial, then an explicit
construction of such
variables $(Y_1, \ldots, Y_{i-1},Y_{i+1}, \ldots,Y_n)$ having the required conditional
distribution, jointly with $(U_1, \ldots U_n)$ is possible.
More details on applications of such constructions to sums of nonlinear functions of Gaussian
and multinomial variables are given in Section \ref{ex}.

\section{Proofs}
\label{pf}
Before proving Theorems \ref{main} and \ref{cor2}, we need the following lemma, the
proof of which can be found in Barbour (1990), or G\"{o}tze (1991).

Let ${\bf Z}$ be a standard $p$-variate normal vector and for $u \ge 0$, define
\[ (T_uh)({\bf w})=E\{ h({\bf w}e^{-u}+\sqrt{1-e^{-2u}}{\bf Z})\}.\]
\begin{lemma}
\label{gotza}
Let $h:R^p\rightarrow R$ have three bounded derivatives. Then
\[ g({\bf w})=-\int_0 ^{\infty} [T_uh({\bf w})-\Phi h]du \]
solves

\[ \mbox{tr}D^2g({\bf w})-{\bf w} \cdot \nabla g({\bf w}) = h({\bf w})-\Phi h,\]
and for any $k^{\rm th}$ partial derivative we have the bound
\[ |\frac{\partial^k}{\prod_{j=1}^k\partial w_{i_j}}g({\bf w})| \leq
\frac{1}{k}\sup_{\bf w}
|\frac{\partial^k}{\prod_{j=1}^k\partial w_{i_j}}h({\bf w})| \le \frac{1}{k}||D^kh||.\]
Further, for any $\blambda \in R^p$ and positive definite $p \times p$ matrix $\Sigma$, $f$ defined by the change of variable
\begin{equation}
\label{solves}
f({\bf w})=g(\Sigma^{-1/2}({\bf w}-\blambda))
\end{equation}
solves
\begin{equation}
\label{DE}
\mbox{tr}\Sigma D^2f({\bf w})-({\bf w}-\blambda) \cdot \nabla f({\bf w}) = h(\Sigma^{-1/2}({\bf w}-\blambda))-\Phi h,
\end{equation}
and hence
\begin{equation}
\label{debound}
|\frac{\partial^k}{\prod_{j=1}^k\partial w_{i_j}}f({\bf w})| \leq
\frac{p^k}{k}||\Sigma^{-1/2}||^k||D^kh||.
\end{equation}
\end{lemma}
{\bf Proof:}\,
One can follow Barbour (1990) to show that $g$ is a solution, and that under the assumptions above, by dominated convergence,
\[ D^kg({\bf w})= -\int_0 ^{\infty} e^{-ku}
E\{ D^kh({\bf w}e^{-u}+\sqrt{1-e^{-2u}}{\bf Z})\}du.\] The Lemma
now follows by straightforward calculations. $\square$

{\bf Proof of Theorem \ref{main}}
Given $h$, let $f$ be the solution of (\ref{DE}) given by (\ref{solves}).
Writing out the expressions in (\ref{DE}) we have
\begin{equation}
\label{first}
 E\{h(\Sigma^{-1/2}({\bf W} - {\blambda}))- \Phi h\}
= E\{\sum_{i=1}^p \sum_{j=1}^p\sigma_{ij}
\frac{\partial^2}{\partial w_i\partial w_j}f({\bf W}) - \sum_{i=1}^p(W_i-\lambda_i)
\frac{\partial}{\partial w_i}f({\bf W}) \}.
\end{equation}
Recall that ${\bf W}^i$ can be characterized by (\ref{rest}):
$$E W_i G({\bf W})=\lambda_i EG({\bf W}^i), $$
holding for all functions $G:R^p \rightarrow R$ for which the expectations exist.
Identity (\ref{rest}) is equivalent to
\begin{equation}
\label{second}
E(W_i-\lambda_i)G({\bf W}) =\lambda_i E[G({\bf W}^i)-G({\bf W})].
\end{equation}
For the coordinate function $G({\bf w})=w_j$ we obtain
\begin{equation}
\label{third}
\sigma_{ij}=\mbox{Cov}(W_i,W_j)=EW_iW_j-\lambda_i \lambda_j=E \lambda_i(W_j^i-W_j);
\end{equation}
when $i=j$ this recovers the one dimensional relation given in Baldi, Rinott, Stein (1989), $E\lambda(W^*-W)=\sigma^2$, where $W^*$ has the $W$-size biased distribution.
Equation (\ref{first}), and (\ref{second}) with $G=\frac{\partial}{\partial w_i}f$, yield
\begin{equation}
\label{fourth}
E\{h(\Sigma^{-1/2}({\bf W} - {\blambda})) - \Phi h\}
= E\{\sum_{i=1}^p \sum_{j=1}^p\sigma_{ij}
\frac{\partial^2}{\partial w_i\partial w_j}f({\bf W}) -
\sum_{i=1}^p\lambda_i
[\frac{\partial}{\partial w_i}f({\bf W}^i)-
\frac{\partial}{\partial w_i}f({\bf W})]\}.
\end{equation}
Taylor expansion of $\frac{\partial}{\partial w_i}f({\bf W}^i)$
centered at ${\bf W}$, with remainder
in integral form, and simple calculations show
that (\ref{fourth}) equals
\begin{eqnarray}
\label{fifth}
\lefteqn{
-E\{\sum_{i=1}^p \sum_{j=1}^p[\lambda_i(W_j^i-W_j)-\sigma_{ij}]
\frac{\partial^2}{\partial w_i\partial w_j}f({\bf W})\}}\\
&&
-E\{\sum_{i=1}^p \sum_{j=1}^p\sum_{k=1}^p\lambda_i
\int_0^1(1-t)\frac{\partial^3}{\partial w_i\partial w_j\partial w_k}
f({\bf W}+t({\bf W}^i-{\bf W}))\,(W_j^i-W_j)(W_k^i-W_k)dt\}. \nonumber
\end{eqnarray}
In the first term, we condition on ${\bf W}$, apply the
Cauchy-Schwarz inequality and use (\ref{third}), and then apply
the bound (\ref{debound}) with $k=2$ to obtain the first term in
(\ref{bounda}). The second term in (\ref{fifth}) gives the second
term in (\ref{bounda}) by applying (\ref{debound}) with $k=3$.
$\square$

{\bf Proof of Theorem \ref{cor2}} Our proof extends and simplifies the proof of Stein~(1986)
in the univariate case.

With ${\bf W} = (W_1,\ldots, W_p)$ where $W_j=\sum_{\alpha \in A_j}X_{\alpha},$
let $W_j^{(\alpha)}=\sum_{\beta \in A_j \cap S_{\alpha}^c}X_{\beta}$,
and ${\bf W}^{(\alpha)}=(W_1^{(\alpha)},\ldots,W_p^{(\alpha)})$.

Let $f$ be the solution of (\ref{DE}) given by (\ref{solves} for a test function $h$.
Writing  (\ref{DE}) while noting that
${\blambda} = {\bf 0}$, and subtracting and adding a term at the end of the expression, we obtain
\begin{eqnarray}
\label{repbounda}
\lefteqn{
 E\{h(\Sigma^{-1/2}{\bf W} ) - \Phi h\}
= E\{\sum_{i=1}^p \sum_{j=1}^p\sigma_{ij}
\frac{\partial^2}{\partial w_i\partial w_j}f({\bf W}) -
\sum_{i=1}^pW_i
\frac{\partial}{\partial w_i}f({\bf W}) \} }\nonumber \\
&& = E\left\{\sum_{i=1}^p \sum_{j=1}^p\sigma_{ij}
\frac{\partial^2}{\partial w_i\partial w_j}f({\bf W}) \right. \nonumber \\
&& \left. -\sum_{i=1}^p\left[ \sum_{\alpha \in A_i} X_{\alpha}[
\frac{\partial}{\partial w_i}f({\bf W})-
\frac{\partial}{\partial w_i}f({\bf W}^{(\alpha)})] + \sum_{\alpha \in A_i}
X_{\alpha}\frac{\partial}{\partial w_i}f( {\bf W}^{(\alpha)})\right] \right\}.
\end{eqnarray}
Taylor expansion of $\frac{\partial}{\partial w_i}f(   {\bf W}^{(\alpha)})$ centered at ${\bf W}$ and some rearrangement shows that the $i^{\rm th}$ summand in the above expression equals
\begin{eqnarray}
\label{taylor}
\lefteqn{
E\{\sum_{j=1}^p
\frac{\partial^2}{\partial w_i\partial w_j}f({\bf W})
[\sigma_{ij} -\sum_{\alpha \in A_i}X_{\alpha}(W_j-W_j^{(\alpha)})]\}
+ E\{\sum_{\alpha \in A_i}X_{\alpha}\frac{\partial}{\partial w_i}f({\bf W}^{(\alpha)})\}} \\
&&
-E\{\sum_{j=1}^p\sum_{k=1}^p\sum_{\alpha \in A_i}X_{\alpha}
\int_0^1(1-t)\frac{\partial^3}{\partial w_i\partial w_j\partial w_k}
f({\bf W}+t({\bf W}^{(\alpha)}-{\bf W}))(W_j-W_j^{(\alpha)})(W_k-W_k^{(\alpha)})dt\}.\nonumber
\end{eqnarray}
Using (\ref{debound}), applying the Cauchy-Schwarz inequality to
the first expectation in (\ref{taylor}), and elementary
calculations on the remaining two terms yield the three terms of
(\ref{secbounda}) respectively. $\square$

\section{Examples}
\label{ex}
\subsection{Sums of nonlinear functions}
Various detailed applications of Theorem \ref{unimain} in the setting of (\ref{fun})
and Lemma \ref{biasfun}, and related references,
are given in Dembo and Rinott (1994). We highlight two problems:
\begin{thm}
\label{gauss}
Let ${\bf U}=(U_1, \ldots, U_n)$ have the multivariate normal
distribution $N({\bf 0}, \Sigma)$, where $\Sigma=\{\rho_{ij}\}$
satisfies $\rho_{ii}=1$
for all $i$, and
$\max_{i \neq j} |\rho_{ij}| \leq r < 1/3$.
Let $W=\sum_{i=1}^n \psi(U_i)$ where
$0 \leq \psi(u) \leq K e^{K |u|^q}$ for some $K>0$ and $q<2$,
and $\psi$ scaled
such that $N\psi=1$ (hence $EW=n$).
Denote $\sigma^2=Var W$. Suppose
$\max_i \sum_{j=1}^n|\rho_{ij}| \leq B < \infty$.
Define $D = \max \{ \|h\|,\|h'\|\}$. Then, for some $C=C(r,K,q)< \infty$,
\begin{equation}
\label{KKn}
|Eh(\frac{W-n}{\sigma})-\Phi h| \leq
 \frac{4 D B (n B C)^{1/2}}{\sigma^2}
+ \frac{D C B^2 n}{\sigma^3}.
\end{equation}
\end{thm}
The construction of $W^*$ utilizes the well known structure of conditional
distributions in the Gaussian case.

An analogous Normal
approximation holds for $W=\sum_{i=1}^n \psi(U_i)$, when
${\bf U}$ is a vector of multinomial variables with equal
(or commensurate) cell probabilities and $\sum_{i=1}^nU_i=kn$ for some
integer $k>0$. In view of Lemma \ref{biasfun} the construction of
the coupling can be done as follows: thinking about $(U_1, \ldots, U_n)$
as counting the distribution of $kn$ balls in $n$ cells,
choose a cell at random, and if cell $i$ is chosen
reset the number of balls in it according to the distribution
$\psi(u)P(U_i \in u)/E\psi(U_i)$. If doing this requires the addition
of balls into cell $i$, these balls are chosen with equal probability
per ball from the other cells. If the resetting requires a reduction
in the number of balls in cell $i$, then a suitable number of balls
is redistributed at random in the remaining cells. Now $W^*$ is the sum
of the function $\psi$ applied to these adjusted cell counts.
This defines $W^*$ on a joint space with $W$, allowing the calculation
of the bound in Theorem \ref{unimain}. This construction generalizes
to any situation where $(U_1, \ldots, U_n)$ have the same distribution
as that of some $n$ iid variables conditioned on their sum.

\subsection{Graph degree counts}
\label{ex1}
Let $K=K_{n,\pi}$ be a random graph on the vertex set
$\{1,\dots ,n\}$, where each pair of vertices has
probability ${\pi}$ of making up an edge, independently of
all other such pairs.  For distinct, fixed $d_i$, $i=1,\ldots,p$, let $W_i$ be the number of vertices of degree $d_i$. Set ${\bf W}=(W_1,\ldots,W_p), \,
E{\bf W}={\blambda}=(\lambda_1,\ldots,\lambda_p)$, \,
$\mbox{Var}{\bf W} = \Sigma = (\sigma_{ij})\,$. For explicit expressions
of ${\blambda}$ and $\Sigma$, see (\ref{cov}) below.
Set
\begin{equation}
\label{beta}
\beta(i)={{n-1} \choose {d_i}}\pi^{d_i} (1-\pi)^{n-1-d_i}, \quad \mbox{and} \quad B=
\bigg[\frac{1}{\min_i\beta(i)
(1-\sum_{i=1}^p\beta(i))}\bigg].
\end{equation}

The theorem below can be extended to the case $n\pi_n\rightarrow
c > 0$  as $n\rightarrow \infty$; for simplicity we assume $0 < \pi=c/(n-1)<1$.

\begin{thm}
\label{graphdeg}
If $\pi=\pi_n=c/(n-1)$, then for any
$h : R^p \rightarrow R$, having bounded mixed partial
derivatives up to order 3,
\begin{equation}
\label{degbounda}
|Eh(\Sigma^{-1/2}({\bf W} - {\blambda})) - \Phi h|  \leq
n^{-1/2}M\bigg\{p^3B||D^2h||(1+c^{5/2})((n-1)/(n-1-c))^{1/2}
+ p^5B^{3/2}||D^3h||(1+c^2)\bigg\}\,,
\end{equation}
where
\begin{equation}
\label{cov}
\lambda_i=n\beta(i), \quad \sigma_{ij}={\mbox Cov}(W_i,W_j) = n\beta(i)\beta(j)
\left[\frac{(d_i-c)(d_j-c)}{c(1-c/(n-1))}-1\right]
+ \delta_{i,j}n\beta(i),
\end{equation}
$\delta_{i,j}$ is 1 if $i=j$, 0 otherwise, and
$M$ is a universal constant. Asymptotic joint normality obviously follows.
\end{thm}

For the case $p=1$, Karo\'nski and Ruci\'nski (1987) proved asymptotic
normality when $\pi n^{(d+1)/d}\rightarrow \infty$ and
$\pi n\rightarrow 0$, or $\pi n \rightarrow \infty$ and $\pi n-\log
n -d \log \log n \rightarrow - \infty$. See also Palka (1984) and Bollob\'{a}s (1985). Asymptotic normality when $\pi n \rightarrow c > 0$, was obtained by Barbour, Karo\'{n}ski and Ruci\'{n}ski (1989). See also Kordecki (1990)
for the case of the one dimensional distribution of the number of vertices of degree zero, for nonsmooth $h$. Numerous univariate results on asymptotic
normality of counts on random graphs, including counts of the type
discussed in Theorems \ref{graphdeg} and 4.3, are given in Janson and Nowicki (1991) and references
therein.

We remark that the calculation of a bound on the conditional variance in Theorem \ref{main}, as well as other terms, is usually involved in nontrivial cases.
The technical details omitted in the following sketch of the proof of Theorem \ref{graphdeg} are available in Goldstein and Rinott (1994).

{\bf Sketch of Proof:} \,
Let $D(v)$ denote the degree of vertex $v$ in $K$ and set $X_{vi}=1$ if
$D(v)=d_i$ and 0 otherwise. We have $W_i=\sum_{v=1}^n X_{vi}, \,\, i= 1,\dots ,p$.
Note that
$D(v) \sim Binomial(n-1, \pi)$ which approaches Poisson$(c)$
as $n \rightarrow \infty$. Note that $\beta(i)$ in (\ref{beta}) equals $EX_{vi} = P(D(v)=d_i)$,
and so the expression for $\lambda_i$ in (\ref{cov}) is immediate. Also, by conditioning on the existence of an edge between vertices $v$ and $u$ and then unconditioning, we can
compute
$EX_{vi}X_{uj}$, and a straightforward calculation
leads to the expression for $\sigma_{ij}$ in (\ref{cov}).
Using the fact that
the maximum absolute value of $||\Sigma^{-1/2}||$ is bounded by the largest eigenvalue of $\Sigma^{-1/2}$, and invoking the Rayleigh-Ritz characterization
of eigenvalues, we obtain with some calculation, $||\Sigma^{-1/2}|| \le N^{-1/2}B^{1/2}$.

The construction required for the application of Theorem
\ref{main} is straightforward.  Fix $i$ and let $V$ be uniformly distributed
on the vertex set $\{1,\dots,n\}$, independent of $K$.
(Note that here $V$ is uniform because
$\gamma_{vi}=EX_{vi}, \, v=1,\ldots, n$ are
all equal for a fixed $i$.) Now $D(V)$ denotes the degree
of the randomly chosen vertex $V$. If $D(V) > d_i$ define $K^i$ to differ from $K$ only in
that $D(V)-d_i$ edges selected uniformly from the $D(V)$ edges
at $V$ are removed from the edge-set. If $D(V) < d_i$ define $K^i$ to be the graph obtained from $K$ by adding $d_i-D(V)$ edges of the form $(V,u)$,
where the vertices $u$ are selected uniformly
from the $n-1-D(V)$ vertices not connected to $V$. If $D(V)=d_i,\
K^i=K$.  Clearly, If $V=v$, then
in the new graph $K^i$ the degree
of $v$ is $d_i$ so that the indicator $X_{vi}$ is
size biased to 1, and the distribution of
$K^i$ is the same as the conditional distribution of $K$
given $X_{vi}=1$.

Define ${\bf W}^i$ to be related to
$K^i$ as ${\bf W}$ is related to $K$, that is,
set
$X^i_{vj}=1$ if in the graph $K^i,$
$D(v)=d_j$, otherwise $0$,
and $W_j^i=\sum_{u=1}^n X^i_{uj}, \,\, j= 1,\dots p$. From the discussion in
Section \ref{secsize} it follows that this
procedure defines ${\bf W}^i$
as in Definition \ref{size} and Theorem \ref{main}.

In order to obtain a tractable bound to the first term on the right hand side of (\ref{bounda}),
we condition on a larger $\sigma$-field, as discussed in Section \ref{secsize}.
Specifically, we use the relation $\mbox{Var}E[\,W_j^i-W_j\, |\, {\bf W}]
\leq \mbox{Var}E[\,W_j^i-W_j\, |\, K],$
and show that
\begin{equation}
\label{convar}
\mbox{Var}E[\,W_j^i-W_j\, |\, K] = O((1+c^3)(1+\frac{d_i^2}{1-c/(n-1)})/n).
\end{equation}

Let $\cal E$ denotes the edge set of $K$ and $|\cdot
|$ cardinality.  Conditioning on $V=v$ and then taking expectation,
recalling $P(V=v)=1/n,$ we obtain
\begin{eqnarray}
\label{one}
\lefteqn{E[\,W_j^i-W_j\, |\, K]}\nonumber \\
&&= {1\over n}\sum\limits_{v:D(v) >d_i}
\left[\bigg|\big\{u:(u,v ) \in {\cal E}, D(u)=
d_j+1\big\}\bigg| - \bigg|\big\{u:(u,v) \in {\cal E},
D(u)=d_j\big\}\bigg|\right]{D(v)-d_i\over
D(v)}\nonumber \\
&&+{1\over n}\sum\limits_{v:D(v) < d_i}
\left[\bigg|\big\{u:(u,v) \not\in {\cal E},\ D(u) =
d_j-1\big\}\bigg| - \bigg|\big\{u:(u,v)\not\in {\cal E},
D(u)=d_j\big\}\bigg|\right]{d_i-D(v)\over n-1-D(v)}\nonumber \\
&& +{1\over n} \bigg|\big\{v:D(v)\not=
d_i\big\}\bigg|\delta_{i,j} -
{1\over n} \bigg|\big\{v:D(v)=
d_j\big\}\bigg|(1-\delta_{i,j})\,.
\end{eqnarray}
To understand the first term, for example, note that that if $V=v$ and $D(v) > d_i$, then $X^i_{uj} - X_{uj}=1$ if $(u,v) \in {\cal E},\ D(u)=d_j+1$, and $(u,v)$ is one of
the $d_i-D(v)$ edges removed at $v$ at random, chosen with
probability $(D(v)-d_i)/ D(v).$

The calculation of a bound on the variance of the expression in (\ref{one})
can be done by computing the covariances between the terms.
They involve conditioning on events to induce independence of terms appearing
as products in the covariances, the use of simple coupling arguments and various moment inequalities.

The bound for the second term on the right hand side of
(\ref{bounda}) is obtained by noting that $|W^i_j-W_j| \le |D(V)
-d_i| +1$ and applying simple calculations related to the Binomial
distribution of $D(V)$. $\square$

\subsection{Graph Vertex Color Matching}
\label{ex2}
Let $G=G_n$ be a fixed regular graph on a vertex set ${\cal V}$ of size $n$, with each vertex $v \in {\cal V}$ of degree $d$. The regularity of $G$ implies the set ${\cal E}$ of edges of $G$ has size $N=nd/2$.
 Let $C=\{1,\ldots,p\}$ be a set
of $p$ colors, and suppose that each vertex $v \in {\cal V}$ is independently assigned color $i$
with probability $\pi_i$. Let $B=
\bigg[\frac{1}{\min_i\pi_i^2
(1-\pi_i)}\bigg]$.

\begin{thm}
Let ${\bf W}=(W_1,\ldots,W_p)$ where $W_i,i=1,\ldots,p$ is the number of edges of $G$ that have both vertices of color $i$. Then    $E{\bf W}=\blambda=N(\pi_1^2,\ldots,\pi_p^2)$ and {\bf Var}$(W)=\Sigma=(\sigma_{ij})$ as given in equations~(\ref{sigcdiag}) and~(\ref{sigcodiag}), and for any
$h : R^p \rightarrow R$, having bounded mixed partial
derivatives up to order 3,
\begin{equation}
\label{colbounda}
|Eh(\Sigma^{-1/2}({\bf W}-\blambda))- \Phi h|  \leq
N^{-1/2}M\bigg\{p^4B||D^2h||d^{3/2}+ p^6B^{3/2}||D^3h||d^2\bigg\}\,,
\end{equation}
where $M$ is a universal constant. Asymptotic joint normality obviously follows.
\end{thm}
{\bf Proof:} \,\, First we will obtain $\Sigma=\mbox{Var}{\bf W}$ in order to bound $||\Sigma^{-1/2}||$. Let $X_{ei}=1$ if edge $e$ has color $i$ on both vertices and 0 otherwise, so $W_i=\sum_{e \in {\cal E}}X_{ei}$ counts the number of edges with both vertices of color $i$. We
have $EX_{ei}=\pi_i^2$, and Var$X_{ei}=\pi_i^2(1-\pi_i^2)$. Given an edge $e$, let $S_e$
denote the set of $2d-1$ edges that share a vertex with $e$, including $e$ itself. For the $2(d-1)$ edges $f \in S_e, f \not = e$, Cov$(X_{ei},X_{fi})=\pi_i^3-\pi_i^4$. For $f \not \in S_e$, this covariance is 0 by independence. Thus, \begin{equation}
\label{sigcdiag}
\sigma_{ii}=\mbox{Var}(W_i)=N\pi_i^2(1-\pi_i^2)+2N(d-1)(\pi_i^3-\pi_i^4).
\end{equation}
For different
colors $i \not = j$ for $f \in S_e,$ we have
Cov$(X_{ei},X_{fj})=-\pi_i^2\pi_j^2$; again, for $f \not \in S_e$ this covariance is 0. Hence,
\begin{equation}
\label{sigcodiag}
\mbox{for $i \not =j$,} \quad
\sigma_{ij}=\mbox{Cov}(W_i,W_j)=-N(2d-1)\pi_i^2\pi_j^2.
\end{equation}

Let $A$ and $H$ be the diagonal matrices with $i^{\rm th}$ diagonal entry $\pi_i^3$, $\pi_i^2-\pi_i^3$ respectively, and let $b$ be a column vector with $i^{\rm th}$ component $\pi_i^2$. Then
$\Sigma=N(2d-1)[A-bb^t]+NH$. In order to show that $\Sigma \succeq NH$, let $D$ be the diagonal
matrix with diagonal entries $\pi_i^{3/2}$, and $g$ the column vector with entries $\pi_i^{1/2}$.
Then $A-bb^t=D(I-gg^t)D$. Since $\Sigma \pi_i=1$, it is easy to see that the smallest eigenvalue
of $I-gg^t$ is 0. Hence, $A-bb^t$ is nonnegative definite and $\Sigma \succeq NH$ is established. It follows that $||\Sigma^{-1/2}|| \le N^{-1/2}B^{1/2}$.

We now apply Theorem \ref{cor2} to the mean zero variables
$X_{ei}-\pi_i^2$. When the square in the first term in the bound
(\ref{secbounda}) is expanded and expectation is taken, most terms
vanish by independence, and because $|S_e|=2d-1$, the number of
summands which do not vanish under the root sign is of the order
$Nd^3$. The second term in (\ref{secbounda}) vanishes, and in the
third term each expectation is of order $d^2$. $\square$

\vspace{3mm}

\baselineskip=12pt

\noindent {\Large\bf References}

\begin{enumerate}

\item[]  Baldi, P. Rinott, Y. (1989). {\sl On normal approximations of distributions in terms of dependency graphs}, Annals of Probability {\bf 17} , 1646-1650.

\item[]  Baldi, P. Rinott, Y.  and Stein, C.  (1989). {\sl A normal approximations for  the number of local maxima of a random function on a graph}, In Probability, Statistics and Mathematics, Papers in Honor of Samuel Karlin.  T. W. Anderson, K.B. Athreya and D. L. Iglehart eds., Academic Press , 59-81.

\item[]  Barbour, A.D. (1990) {\sl Stein's method for diffusion approximations},
Probab. Th. Rel. Fields {\bf 84} 297-322.

\item[] Barbour, A. D., Karo\'nski, M. and Ruci\'nski, A. (1989). {\sl A central limit theorem for decomposable random variables with applications to random graphs}, J. Combinatorial Theory B {\bf 47}, 125-145.

\item[] Bollob\'as, B.  (1985) {\sl Random Graphs}.  Academic Press, 1985.

\item[] Brewer, K. and Hanif, M. (1983) {\sl Sampling with unequal probabilities }, Lecture Notes in Statistics, vol. 15. Springer-Verlag, New York.

\item[] Cochran, W. (1977) {\sl Sampling Techniques} John Wiley \& Sons, New York.

\item[]  Dall'Aglio, S, Kotz, S., and Salinetti, G. (Eds.) (1991)
{\sl Advances in Probability Distributions With Given Marginals}.
Kluwer Academic Publishers, Dordrecht, Boston, London.

\item[]  Dembo, A, and Rinott, Y. (1994). {\sl Some examples of Normal approximations by Stein's method}. To appear in IMA Conference Proceedings,
Aldous, D. and Pemantle, R. Eds.

\item[] Goldstein, L., and Rinott, Y. (1994) {\sl On multivariate normal
approximations by Stein's method and size bias couplings: Technical
Report}.

\item[] G\"otze, F. (1991) {\sl On the rate of convergence in the multivariate CLT}. Annals of Probability, {\bf 19}, 724-739.

\item[]  Horn, R. A., and Johnson, C. A.  (1985) {\sl Matrix Analysis}. Cambridge University Press 1985.

\item[] Janson, S. and Nowicki, K. (1991) {\sl The asymptotic distributions
of generalized U-statistics with applications to random graphs}.
Probability Theory and Related Fields {\bf 90}, 341-375.

\item[] Karo\'nski, M.  and Ruci\'nski A. (1987), {\sl Poisson
convergence of semi-induced properties of random graphs}. Math. Proc. Comb. Phil. Soc. {\bf 101}
291-300.

\item [ ] Kordecki, W. (1990) {\sl Normal approximation and isolated vertices in random graphs}. In Random Graphs 1987,  M. Karo\'nski, J. Jaworski and A. Ruci\'nski eds. John Wiley \& Sons, New York.

\item [ ] Luk, H. M. (1994) {\sl Stein's method for the Gamma distribution and related statistical applications}. PhD dissertation, USC.

\item [] Palka, Z. (1984)  {\sl On the number of vertices of a given degree in a random graph}.  J. Graph Theory {\bf 8}, 167-170.

\item [ ] Reinert, G. (1994) {\sl A weak law of large numbers for
empirical measures via Stein's method, and applications}. PhD dissertation.

\item[] Rinott, Y. (1994) {\sl On normal approximation rates for certain sums of dependent random variables}. To appear in Journal of applied and Comp. Math.

\item[] Rinott, Y. and Rotar, V. (1994) {\sl  A multivariate CLT for local dependence with $n^{-1/2}logn$ rate.}

\item[] Stein, C. (1972) {\sl A bound for the error in the normal approximation to
the distribution of a sum of dependent random variables.} Proc. Sixth Berkeley Symp. Math. Statist. Probab. {\bf 2} 583-602, Univ. California Press, Berkeley.

\item[] Stein, C. (1986)  {\sl Approximate Computation of Expectations}.   IMS, Hayward, Calif., 1986.

\item[] Stein, C. (1992) {\sl A way of using auxiliary randomization} Probability Theory, pp. 159-180. Walter de Gruyter \& Co., Berlin - New York.

\end{enumerate}

\end{document}